\documentclass[conference, 10pt]{IEEEtran}

\usepackage{cite}
\usepackage{url}
\usepackage{graphicx}
\usepackage{color}
\usepackage{placeins}
\usepackage{float}
\usepackage{tabularx,colortbl}
\usepackage{ifthen}

\usepackage{lipsum}
\usepackage{subfigure}

\usepackage{calc,ifthen,intcalc}

\usepackage{tikz}
	\usetikzlibrary{arrows,shapes,positioning,calc,fadings,patterns,decorations}
	\usetikzlibrary{decorations.pathreplacing,decorations.pathmorphing}			
	\usetikzlibrary{arrows,arrows.meta}
	\usetikzlibrary{shapes.geometric}

\newcommand{\latexConstant}{28.45}
\newcommand{\xOffset}{0.9}
\newcommand{\xComplement}{0.1}
\newcommand{\xScale}{0.4}
\newcommand{\yScale}{\xScale*\xOffset}

\makeatletter

\newcounter{author}
\renewcommand{\author}[2][]{
   \stepcounter{author}
   \@namedef{author@\theauthor}{#2}
   \@namedef{authorlabel@\theauthor}{#1}
}

\newcounter{address}
\newcommand{\address}[2][]{
   \stepcounter{address}
   \@namedef{address@\theaddress}{#2}
   \@namedef{addresslabel@\theaddress}{#1}
}

\newcommand{\alsep}{and}

\def\newmaketitle{\par%
  \begingroup%
  \normalfont%
  \def\thefootnote{}
  \def\footnotemark{}
  \let\@makefnmark\relax
  \footnotesize
  \footnotesep 0.7\baselineskip
  \normalsize%
  \twocolumn[\thenewmaketitle\@IEEEaftertitletext]%
  \if@IEEEusingpubid
     \enlargethispage{-\@IEEEpubidpullup}%
  \fi
  \endgroup
  \setcounter{footnote}{0}\let\maketitle\relax\let\@maketitle\relax
  \gdef\@thanks{}%
  \let\thanks\relax}

\def\thenewmaketitle{
  \newpage
  \begin{center}%
    \vskip0.2em{\Huge\@IEEEcompsoconly{\sffamily}\@IEEEcompsocconfonly{\normalfont\normalsize\vskip 2\@IEEEnormalsizeunitybaselineskip
   \bfseries\large}\@title\par}\vskip1.0em\par%
    \vspace{1ex}
    \newcounter{c@author}
    \newcounter{c@tmp}
    \ifthenelse{\value{author}=2}{%
      \newcommand{\liand}{ and }}{%
      \newcommand{\liand}{, and }}
    \ifthenelse{\value{address}<2}{%
      \@nameuse{author@1}%
      \stepcounter{c@author}%
      \whiledo{\value{c@author}<\value{author}}{%
        \setcounter{c@tmp}{\value{author}}%
        \addtocounter{c@tmp}{-\value{c@author}}%
        \ifthenelse{\value{c@tmp}=1}{%
          \renewcommand{\alsep}{\liand}}{\renewcommand{\alsep}{, }}%
        \stepcounter{c@author}\alsep \@nameuse{author@\thec@author}}\\%
    }
    {
      \@nameuse{author@1}${}^{(\ref{\@nameuse{authorlabel@1}})}$%
      \stepcounter{c@author}%
      \whiledo{\value{c@author}<\value{author}}{%
      \setcounter{c@tmp}{\value{author}}%
      \addtocounter{c@tmp}{-\value{c@author}}%
      \ifthenelse{\value{c@tmp}=1}{%
        \renewcommand{\alsep}{\liand}}{\renewcommand{\alsep}{, }}%
      \stepcounter{c@author}\alsep \@nameuse{author@\thec@author}%
        ${}^{(\ref{\@nameuse{authorlabel@\thec@author}})}$%
      }
    }
    \vspace{0.2ex}

    \ifthenelse{\value{address}>0}{%
      \ifthenelse{\value{address}=1}{
        {\@nameuse{address@1}}
      }
      {
        \newcounter{c@address}

        \begin{center}
        \whiledo{\value{c@address}<\value{address}}
        {
          \refstepcounter{c@address}
            ${}^{(\thec@address)}$\,%
              \label{\@nameuse{addresslabel@\thec@address}}%
              \@nameuse{address@\thec@address}\\ %
        }
        \end{center}
      } 
    }
    {
      \relax
    }
  \end{center}
}

\usepackage{booktabs,bm,amsmath,amssymb}

\newcommand{\bs}[1]{\bm{\mathrm{#1}}}
\newcommand{\abs}[1]{\left|#1\vphantom{f}\right|}

\renewcommand{\d}{\ensuremath{\partial}}
\newcommand{\e}{\ensuremath{\mathrm{e}}}

\newcommand\figref[1]{Fig.~\ref{#1}}

\makeatother

\title{Combining Simulation and Far-Field Approximation for RIS Optimization}
\title{Efficient Far-Field Approximation for RIS Optimization: Streamlining Analysis with a Single Simulation}
\title{RIS Optimization with Efficient Far-Field Approximation Requiring Only a Single Simulation}
\title{Reducing Simulation Effort using an Efficient Far-Field Approximation for RIS Optimization }
\title{Reducing Simulation Effort for RIS Optimization using an Efficient Far-Field Approximation}

\author[org1]{Hans-Dieter Lang}
\author[org1]{Michel A. Nyffenegger}
\author[org1]{Heinz Mathis}
\author[org2]{Xingqi Zhang}

\address[org1]{OST – Eastern Switzerland University of Applied Sciences, Rapperswil, SG, Switzerland. Email: hansdieter.lang@ost.ch}
\address[org2]{School of Electrical and Electronic Engineering, University College Dublin, Dublin, Ireland.}

\begin{document}
\bstctlcite{IEEEexample:BSTcontrol}

\newmaketitle

\begin{abstract}
Optimization of Reconfigurable Intelligent Surfaces (RIS) via a previously introduced method is effective, but time-consuming, because multiport impedance or scatter matrices are required for each transmitter and receiver position, which generally must be obtained through full-wave simulation. Herein, a simple and efficient far-field approximation is introduced, to extrapolate scatter matrices for arbitrary receiver and transmitter positions from only a single simulation while still maintaining high accuracy suitable for optimization purposes. This is demonstrated through comparisons of the optimized capacitance values and further supported by empirical measurements.
\end{abstract}

\vspace*{4pt}
\section{Introduction}

Reconfigurable Intelligent Surfaces (RISs) are recognized as key components in the evolution of wireless communication technologies, moving from "beyond 5G" to 6G and subsequent generations~\cite{renzo20,pan21}. 
They intelligently redirect wireless signals otherwise lost or “spilled” back into the system, thereby enhancing both coverage and throughput of wireless connections. 

To assess the performance of an RIS in a particular scenario against the limit of what is possible, an effective and versatile optimization approach has been introduced \cite{lang23}. It is capable of determining the optimal reactive loading components for $N$ ports of reconfigurable antenna elements by maximizing power transfer from a specific transmitter (Tx), the source of the impinging signal,  
to a designated receiver (Rx).  


This technique relies on the impedance or scatter matrix of the entire system, necessitating a multiport simulation for each desired Tx and Rx position. Even with fast simulation methods such as the Multilevel Fast Multipole Method (MLFMM) in ANSYS HFSS, these simulations can become time-intensive. To address this, a simple far-field approximation is presented that allows the extrapolation of scatter matrices for all positions from just one single multiport simulation. The validity of this approach is confirmed through comparisons between simulations, approximations, and measurements.

\vspace*{4pt}
\section{Far-Field Approximation}

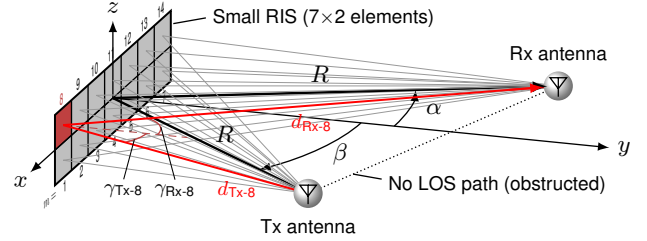
\begin{figure}[!htb]\centering
	\vspace*{2mm}
	\begin{tikzpicture}[scale=1.1,>=latex,line width=0.5pt]\sffamily\small
		
		\newcommand{\shiftX}{0*\latexConstant}
\begin{scope}[yshift=\shiftX*-\xComplement,xshift=\shiftX]
\begin{scope}[every node/.append style={yslant=\xOffset,xscale= \xScale},yslant=\xOffset,xscale= \xScale,xshift=0,yscale=0.85]

	\fill[fill=lightgray] (-1.75,-0.5)rectangle(1.75,0.5);
	\fill[red!50!gray] (-1.75,0)rectangle(-1.25,0.5);
	\draw[line width=0.75pt] (-1.75,-0.5)rectangle(1.75,0.5);
	\draw (-1.75,0)--(1.75,0);
	\foreach \x in {-1.25,-0.75,...,1.25}{
	\draw (\x,-0.5)--(\x,0.5);
	}
	
	\foreach \n in {1,2,...,7}{
		\node[scale=0.66] at(\n*0.5-1.95,-0.64){\n};
	}
	\def\n{8}
		\node[scale=0.66,red!50!gray] at({(\n-7)*0.5-2.05},0.64){\n};
	\foreach \n in {9,10,...,14}{
		\node[scale=0.66] at({(\n-7)*0.5-2.05},0.64){\n};
	}
	
	\node[scale=0.66] at(-1.9,-0.64){$m=$};
	
	\node[coordinate] (normal) at(1.5,0.25){};
	\node[coordinate] (normal2) at(-1.5,0.25){};
	
	\sffamily
	
	\draw[->] (0,0)--(-2.5,0)coordinate[below left=5pt](xtemp);
	\draw[->] (0,0)--(0,1.1)coordinate[above=5pt](ztemp);

\end{scope}
\end{scope}

\renewcommand{\shiftX}{0*\latexConstant}

\begin{scope}[xshift=-\xScale*\shiftX,yshift=-\yScale*\shiftX]
\begin{scope}[every node/.append style={yslant=\xOffset-1},yslant=\xOffset-1]	

	
		\draw[->] (0,0)--(6,0)node[right]{$y$};
	\node at(xtemp){$x$};
	\node at(ztemp){$z$};

\end{scope}
\end{scope}

\def\aalpha{35}
\def\bbeta{-20}

\begin{scope}[yshift=0*\latexConstant,scale=-2]
\begin{scope}[every node/.append style={yslant=\xOffset,xslant=-1,xscale= \xScale},yslant=\xOffset,xslant=-1,xscale= \xScale]

\def\delta{0.1758}
		\def\T{2.5}
		\def\R{2}
			
		
		\path[miter limit=1] (90-\bbeta:\T)coordinate(bbeta)--(0,0)--(90-\aalpha:\R)coordinate(aalpha);

\end{scope}
\end{scope}

\begin{scope}[yshift=\shiftX*-\xComplement,xshift=\shiftX]
\begin{scope}[every node/.append style={yslant=\xOffset,xscale= \xScale},yslant=\xOffset,xscale= \xScale,xshift=0]

			\def\n{-1}
		\foreach \m in {-2.5,-1.5,...,3.5}{
			\draw[gray!80!white,miter limit=1,line width=0.33pt] (bbeta)--(-0.25+0.5*\m,{(0.25+\n*0.5)*0.85})--(aalpha);
		}
		
\end{scope}
\end{scope}

\begin{scope}[yshift=0*\latexConstant,scale=-2]
\begin{scope}[every node/.append style={yslant=\xOffset,xslant=-1,xscale= \xScale},yslant=\xOffset,xslant=-1,xscale= \xScale]

\def\delta{0.1758}
		\def\T{2.5}
		\def\R{2}
			
		
		\draw[line width=0.9pt,miter limit=1] (90-\bbeta:\T)coordinate(bbeta)--(0,0)--(90-\aalpha:\R)coordinate(aalpha);
		
		\draw[->] (90:\T-0.8)arc(90:90-\bbeta:\T-0.8);
		\coordinate (bet) at(90-0.5*\bbeta:\T-0.66){};
		\draw[->] (90:\R-0.5)arc(90:90-\aalpha:\R-0.5);
		\coordinate (alp) at(90-0.5*\aalpha:\R-0.35){};

\end{scope}
\end{scope}

\begin{scope}[yshift=\shiftX*-\xComplement,xshift=\shiftX]
\begin{scope}[every node/.append style={yslant=\xOffset,xscale= \xScale},yslant=\xOffset,xscale= \xScale,xshift=0]

		
			\def\n{0}
		\foreach \m in {-1.5,-0.5,...,3.5}{
			\draw[gray!80!white,miter limit=1,line width=0.33pt] (bbeta)--(-0.25+0.5*\m,{(0.25+\n*0.5)*0.85})--(aalpha);
		}
		

\end{scope}
\end{scope}

\renewcommand{\shiftX}{0*\latexConstant}

\begin{scope}[xshift=-\xScale*\shiftX,yshift=-\yScale*\shiftX]
\begin{scope}[every node/.append style={yslant=\xOffset-1},yslant=\xOffset-1]	
	
	\draw[densely dashed,line width=0.4pt,red!40!gray] (normal2)--coordinate[pos=0.67](tempp)coordinate[pos=0.76](temppp)+(1.5,0);

\end{scope}
\end{scope}

\begin{scope}[yshift=0*\latexConstant,scale=-2]
\begin{scope}[every node/.append style={yslant=\xOffset,xslant=-1,xscale= \xScale},yslant=\xOffset,xslant=-1,xscale= \xScale]
		
		\draw[line width=0.5pt,red!50!gray] (tempp)arc(90:49:0.25);
		\draw[line width=0.5pt,red!50!gray] (temppp)arc(90:155:0.25);

\end{scope}
\end{scope}

\begin{scope}[yshift=\shiftX*-\xComplement,xshift=\shiftX]
\begin{scope}[every node/.append style={yslant=\xOffset,xscale= \xScale},yslant=\xOffset,xscale= \xScale,xshift=0]

		\def\n{0}
		\def\m{-2.5}
			\draw[miter limit=1,red,line width=0.75pt] (bbeta)--coordinate[pos=0.1](aaalpha)(-0.25+0.5*\m,{(0.25+\n*0.5)*0.85})--(aalpha);

\end{scope}
\end{scope}
	
	\node at(alp){$\beta$};
	\node at(bet){$\alpha$};
	
	\node at(2.5,0.29){$R$};
	\node at(1.35,-0.43){$R$};
	\node[red,scale=0.85] at(1.5,-1.1){$d_\text{Tx-8}$};
	\node[red,scale=0.85] at(2.4,-0.3){$d_\text{Rx-8}$};
	
	\draw[densely dotted] (bbeta)--coordinate[pos=0.82](temp)(aalpha);
	\draw (temp)--++(-10:0.35)+(0,-2pt)node[right,scale=0.8]{No LOS path (obstructed)};
	
	\node[circle,ball color=white,minimum size=11pt] at(bbeta){};
	\node[circle,ball color=white,minimum size=11pt] at(aalpha){};
	
	\draw (bbeta)+(0,-0.13)--++(0,0.1)--+(-0.1,0)--+(0,-0.133)--+(0.1,0)--+(0,0);
	\draw (aalpha)+(0,-0.13)--++(0,0.1)--+(-0.1,0)--+(0,-0.133)--+(0.1,0)--+(0,0);
	\node[scale=0.8,above=7pt] at(bbeta){Rx antenna};
	\node[scale=0.8,below=7pt] at(aalpha){Tx antenna};
	
	\draw[line width=0.4pt] (0.26,-0.5)--+(-105:0.5)node[below=-1pt,scale=0.85]{$\gamma_\text{Tx-8}$};
	\draw[line width=0.4pt] (0.58,-0.35)--+(-75:0.66)node[below=-1pt,scale=0.85]{$\gamma_\text{Rx-8}$};

	\path (bbeta)--coordinate[pos=0.025](tempa)coordinate[pos=0.05](tempb)(0,0);
	\draw[<-,line width=0.9pt] (tempa)--(tempb);
	
	\path (aaalpha)--coordinate[pos=0.7](aaaalpha)(bbeta);
	\draw[red,->,line width=0.75pt] (aaalpha)--(aaaalpha);
	
	\draw (0.73,0.9)--+(10:0.4)node[right,scale=0.8] {Small RIS (7$\times$2 elements)};

	\end{tikzpicture}
	\vspace*{-1.5mm}
	\caption{Geometric illustration of the setup: The RIS facilitates an NLOS communication channel between the transmitter and the receiver antennas, positioned at angles $\beta$ and $\alpha$ w.r.t. the surface normal, respectively. The red path highlights the coupling of the Tx via RIS (at element $m=8$) to the Rx.}
	\label{fig:overall_setup}
	\vspace*{-1mm}
\end{figure}

\figref{fig:overall_setup} illustrates the typical RIS-aided NLOS (non-line-of-sight) wireless link setup. Optimizing this link by finding optimal loading reactances for each tunable antenna element on the RIS using the maximum power transfer framework \cite{lang23}, requires impedance or scatter matrices of the entire system, generally obtained from full-wave simulation. Different incident/Tx angles $\beta$ and reflection/Rx angles $\alpha$ result in distinct matrices --- however, most matrix entries remain the same, only the Tx-RIS and Rx-RIS coupling entries change. 


For these purposes and sufficient distance between the $m$th RIS antenna element and the Tx, the respective coupling entries of the scatter matrix of the entire system can be estimated as follows:
\begin{equation}
	S_{\text{Tx-}m}\approx \sqrt{1-\abs{S_{mm}}^2}\,\frac{\sqrt{G_\text{Tx}G_m(\gamma_{\text{Tx-}m})}}{4\pi d_{\text{Tx-}m}/\lambda}\,\e^{-j\,2\pi d_{\text{Tx-}m}/\lambda}\,,
\end{equation}
where  $G_\text{Tx}$ and $G_m$ are the gains of the transmit and $m$th RIS antenna element, respectively, and $S_{mm}$ is the scatter coefficient of that element. Essentially, this treats the NLOS link connecting the Tx antenna via $N$ RIS elements to the Rx antenna as $N$ lossy two-ports \cite{collin01} with delay. The first square root  accounts for impedance mismatch on the RIS side; the Tx and Rx antennas are  assumed to be perfectly matched. 

The gain pattern $G_m(\phi)$ of the $m$th RIS antenna element is extracted from the single full-wave simulation when exciting that element separately. The azimuth angle $\gamma_{\text{Tx-}m}$ to the Tx for that $m$th RIS element is approximated from the geometric relation
$\sin\gamma_{\text{Tx-}m}=\left(R\,\sin\beta-x_m\right)/d_{\text{Tx-}m}$\,. Lastly, the distance between the Tx antenna and the $m$th RIS antenna element is calculated by
	$d_{tm}=\sqrt{R^2+x_m^2+z_m^2-2\,x_mR\sin\beta}$\,, where $x_m$ and $z_m$ are the coordinates of the element's position.

\newcommand\npt{\hspace*{-1pt}}
Note that, as $R\npt\rightarrow\npt\infty$, the parallel-ray approximation results, with $d_{\text{Tx-}m}\npt\rightarrow\npt R$ and $\gamma_{\text{Tx-}m}\npt\rightarrow\npt\beta$. However, for later comparison to practical measurements at finite/shorter distances, these formulations prove to be more useful and accurate.

Applying the aforementioned procedure also to the coupling between the RIS and the Rx antenna, where the angle $-\alpha$ replaces $\beta$, leading to $d_{\text{Rx-}m}$ and $\gamma_{\text{Rx-}m}$, the entire $(N+2)\times(N+2)$ 
S-matrix can be systematically constructed as follows:
\renewcommand\npt{\hspace*{0pt}}
\begin{equation}
\npt\bs{S}
		\approx
		\begin{bmatrix}
			S_\text{Tx}\approx0 & \npt\bs{s}_\text{Tx-RIS}^T(\beta)\npt & S_\text{Tx-Rx}\approx 0\\
			\bs{s}_\text{Tx-RIS}(\beta)\npt & \bs{S}_\text{RIS} & \npt\bs{s}_{\text{Rx-RIS}}(\alpha)\\
			S_\text{Tx-Rx}\approx 0 & \npt\bs{s}_\text{Rx-RIS}^T(\alpha)\npt & S_\text{Rx}\approx 0
		\end{bmatrix}
		\,.
		\end{equation}
Hence, only a single full-wave simulation for $\bs{S}_\text{RIS}$ is required, from which the entire $\bs{S}$ can be approximated for all $\alpha$ and $\beta$.

To obtain the desired NLOS scenario with $S_\text{Tx-Rx}\approx0$, well-matched ($S_\text{Tx}\approx S_\text{Rx}\leq -16$\,dB) horn antennas with $G_\text{Tx}\approx G_\text{Rx}\approx 11$\,dB gain   are used as Tx and Rx antennas both in simulation and for measurements. Furthermore, time gating is employed during measurements to reduce influences of unwanted reflections of the measurement environment \cite{phumvijit17}. 

\vspace*{4pt}
\section{Comparison}
For comparison purposes, a small $7\times2$-element RIS for $f=3.55$\,GHz is realized, shown in \figref{fig:capvalues}(a). Both the Tx and Rx antennas are located at $R=2$\,m distance. The RIS is optimized for an Rx at $\alpha=0$°, when the Tx is at $\beta=30$°.

\subsection{Optimized Capacitance Values}
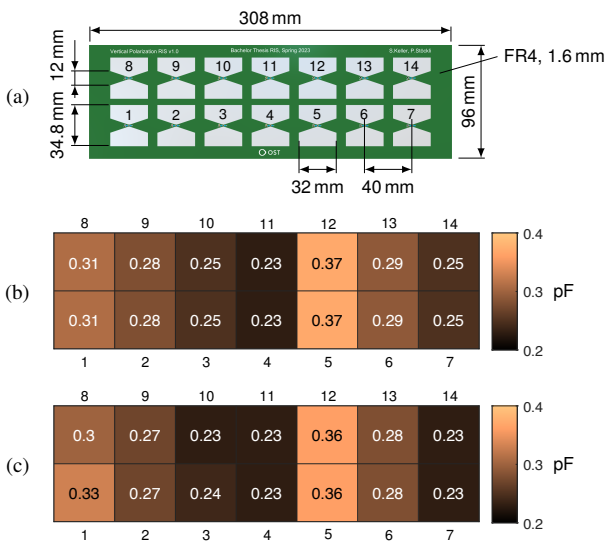
\begin{figure}[!htb]\centering
	\vspace*{-0.5mm}
	\begin{tikzpicture}[>=latex]\small\sffamily
	
	\begin{scope}[yshift=2.55cm,xshift=3cm,line width=0.3pt]
		\node at(-0,0){\includegraphics[width=4.82cm]{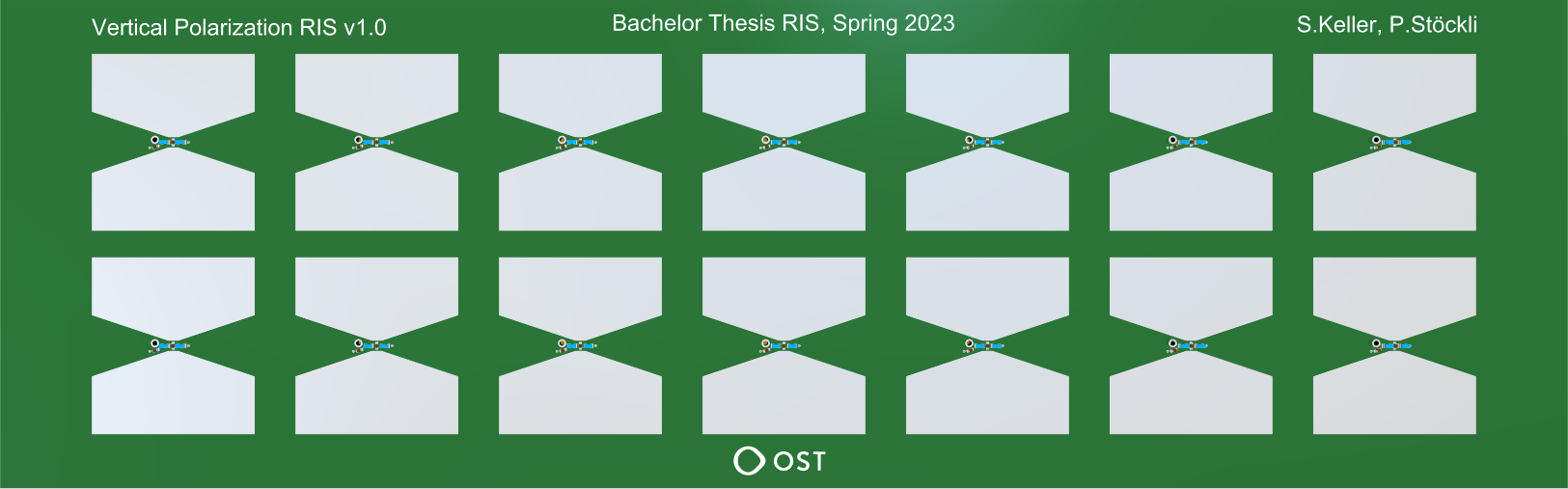}};
		\draw[<->] (-2.4,0.95)--node[above,scale=0.75]{308\,mm}(2.4,0.95);
		\draw (-2.4,1.0)--(-2.4,0.8) (2.4,1.0)--(2.4,0.8);
		\draw[<->] (2.8,-0.75)--node[above,scale=0.75,rotate=90]{96\,mm}(2.8,0.75);
		\draw (2.5,-0.75)--(2.85,-0.75) (2.5,0.75)--(2.85,0.75);
		
		\draw[<->] (1.25,-0.9)--node[below=2pt,scale=0.75]{40\,mm}(1.88,-0.9);
		\draw (1.25,-0.95)--+(0,0.72) (1.88,-0.95)--+(0,0.72);

		\draw[<->] (0.38,-0.9)--node[below=2pt,scale=0.75]{32\,mm}(0.875,-0.9);
		\draw (0.38,-0.95)--+(0,0.43) (0.875,-0.95)--+(0,0.43);
		
		\draw[<->] (-2.6,-0.582)--coordinate(temp)(-2.6,-0.04);
		\path (temp)+(-2pt,0)node[above,scale=0.75,rotate=90]{34.8\,mm};
		\draw (-2.15,-0.582)--(-2.65,-0.582) (-2.15,-0.04)--(-2.65,-0.04);
		
			\draw (-2.6,0.22)--(-2.6,0.41);
			\draw[<-] (-2.6,0.41)--coordinate(temp)+(0,0.41);
			\draw[->] (-2.6,0.05)--(-2.6,0.22);
		\path (temp)+(-2pt,0)node[above,scale=0.75,rotate=90]{12\,mm};
		\draw (-2.15,0.22)--(-2.65,0.22) (-2.15,0.41)--(-2.65,0.41);
		
		\foreach \n in {1,...,7}{
			\node[scale=0.6] at(\n*0.625-0.625-1.88,-0.17){\n};
		}
		\foreach \n in {8,...,14}{
			\node[scale=0.6] at(\n*0.625-0.625*8-1.88,0.47){\n};
		}
		
		\draw[-] (2.25,0.38)--+(16:0.8)node[right,scale=0.75]{FR4, 1.6\,mm};
		
		\end{scope}

	\def\z{2.3}
	\def\x{0.1}
	{\footnotesize\normalfont
	\node[anchor=east,scale=1] at(-0.1,2.6){(a)};
	\node[anchor=east,scale=1] at(-0.1,0){(b)};
	\node[anchor=east,scale=1] at(-0.1,-\z){(c)};
	}
	\node[anchor=west,inner sep=0.1pt] at(\x,0){\includegraphics[scale=0.6,clip,trim={17mm 4mm 6mm 1.5mm}]{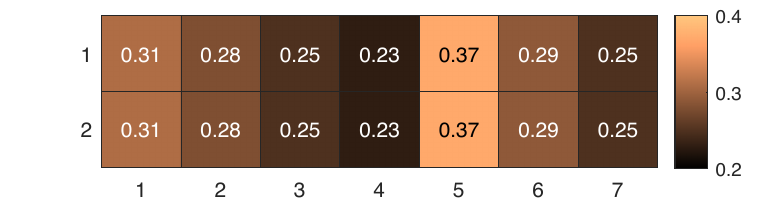}};
		\foreach \n in {1,2,...,7}{
			\node[scale=0.6] at({0.81*(\n-1)+0.53},-0.91){\n};}
		\foreach \n in {8,9,...,14}{
			\node[scale=0.6] at({0.81*(\n-8)+0.53},0.93){\n};}			
		\node[scale=0.75] at(6.9,0+0.01){pF};
	
	\node[anchor=west,inner sep=0.1pt] at(\x,-\z){\includegraphics[scale=0.6,clip,trim={17mm 4mm 6mm 1.5mm}]{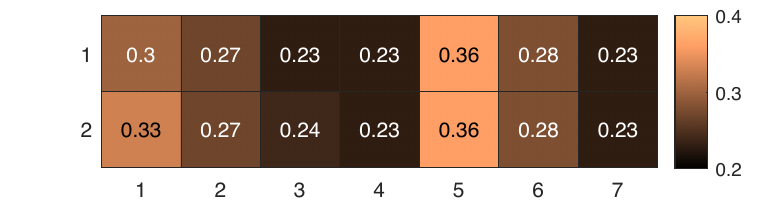}};
		\foreach \n in {1,2,...,7}{
			\node[scale=0.6] at({0.81*(\n-1)+0.53},-\z-0.91){\n};}
		\foreach \n in {8,9,...,14}{
			\node[scale=0.6] at({0.81*(\n-8)+0.53},-\z+0.93){\n};}
		\node[scale=0.75] at(6.9,-\z+0.01){pF};
	\end{tikzpicture}
	\vspace*{-1mm}
	\caption{The RIS (a) and a comparison of the resulting capacitance values for (b) the approximation and (c) the full simulation (using FEBI with MLFMM in ANSYS HFSS), for $\alpha=0$°, $\beta=30$°.}
	\label{fig:capvalues}
	\vspace*{-1mm}
\end{figure}

 The optimized capacitance values are illustrated for the approximation in \figref{fig:capvalues}(b) and for the full simulation in (c). As can be seen, the obtained values are very close.

\subsection{Measurements}
The capacitance values are realized using reverse-biased Skyworks SMV2201-040LF varactor diodes \cite{SMV2201} on the RIS. The achievable range of capacitance is limited to approx. 0.23 to 2.1\,pF, adequate for realizing the values obtained in \figref{fig:capvalues}.

The obtained bistatic radar cross sections (BRCS) [3] show strong agreement across the approximation, full simulation, and measurements, confirming the practical viability of the approximation and the realizability of the optimization results.

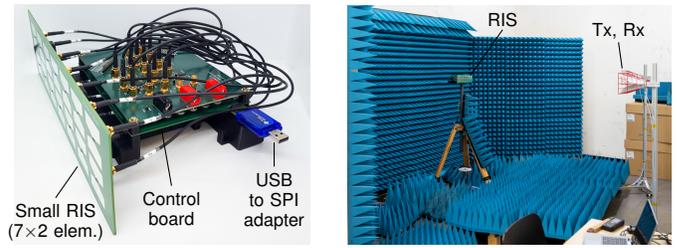
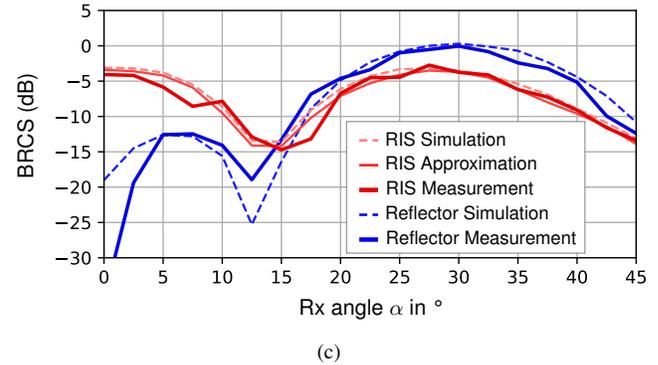
\begin{figure}[!htb]\centering
	\subfigure[]{\begin{tikzpicture}[scale=0.5335]\sffamily\small
		\node[inner sep=0pt] at(0,0){\includegraphics[width=3.93cm,clip,trim={0 1mm 0 12mm}]{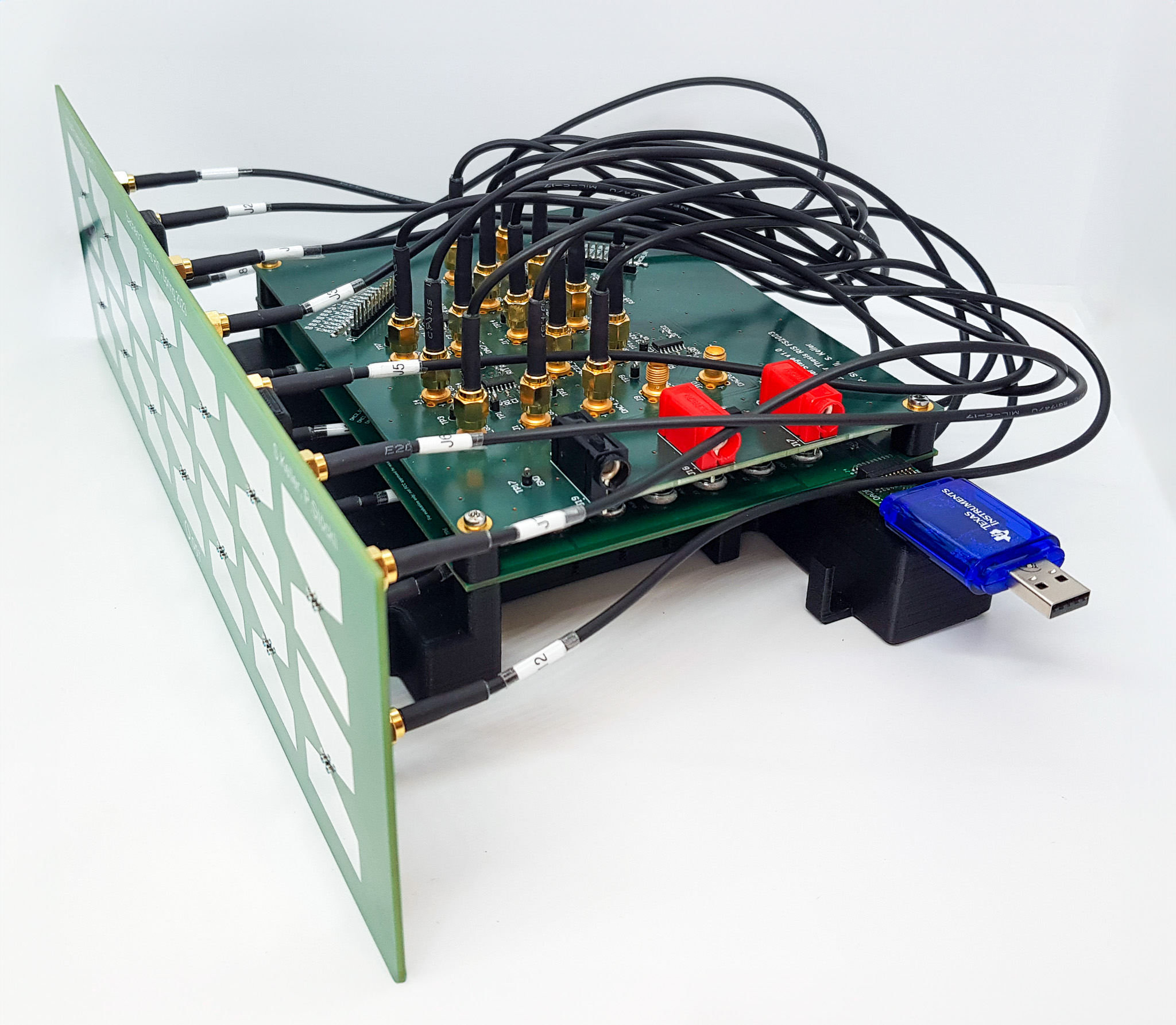}};
		\draw (-2.1,-1.05)--+(-125:0.95)node[text width=1.8cm,align=center,scale=0.75,below=2pt,inner sep=0pt]{Small RIS\\ (7$\times$2 elem.)};
		\draw (0,-0.15)--+(-80:1.4)node[below,text width=1.6cm,align=center,scale=0.75]{Control board};
		\draw (2.8,-0.45)--+(-90:0.6)node[below,text width=1.1cm,align=center,scale=0.75]{USB to SPI adapter};
	\end{tikzpicture}}
	\hfill
	\subfigure[]{\begin{tikzpicture}\sffamily\small
		\node[inner sep=0pt]{\includegraphics[width=4.3cm,clip,trim={26mm 2mm 9mm 16mm}]{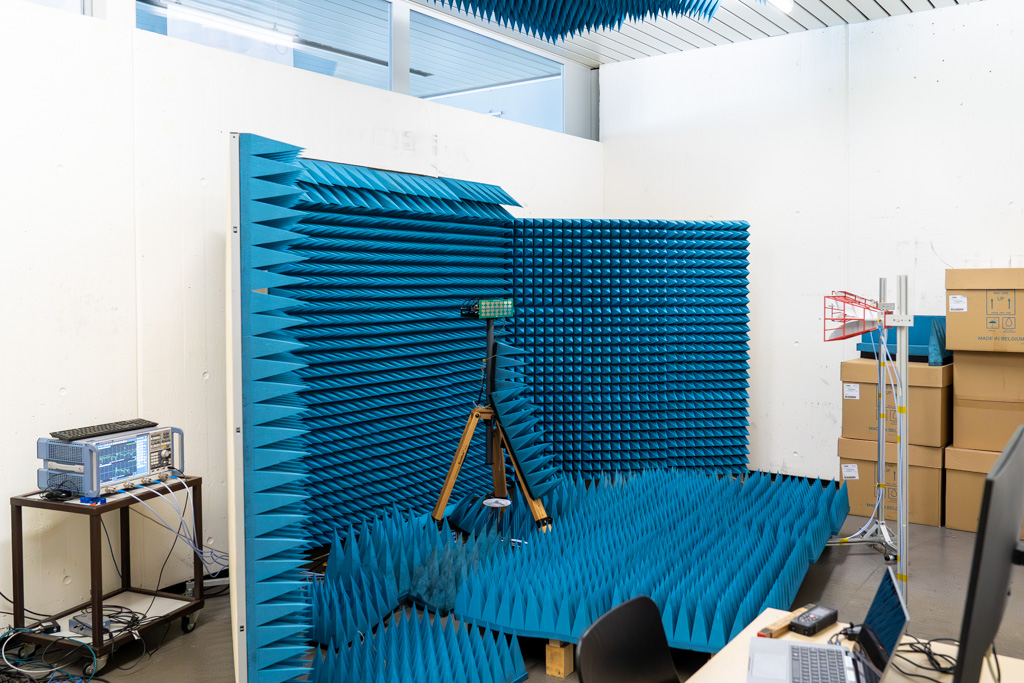}};
		\draw (1.5,0.75)--+(98:0.35)node[above=-1pt,text width=1.1cm,align=center,scale=0.75]{Tx, Rx};
		\draw (-0.52,0.65)--+(55:0.75)node[above=-1pt,text width=1.1cm,align=center,scale=0.75]{RIS};
		\end{tikzpicture}}\\[1mm]
	\subfigure[]{\begin{tikzpicture}\sffamily\small
		\node[inner sep=0pt] at(0,0){\includegraphics[scale=0.66,clip,trim={2mm 4mm 2mm 2mm}]{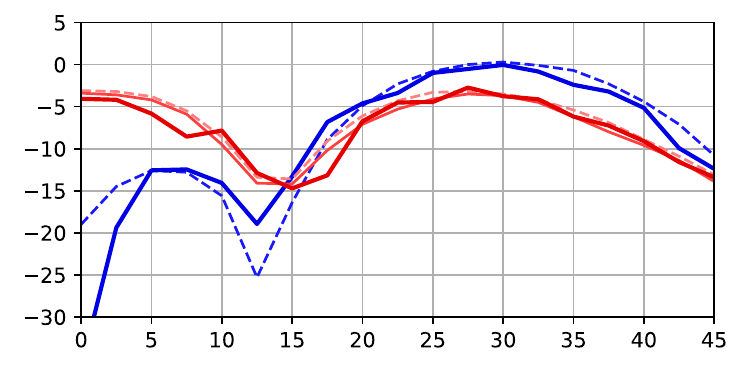}};
		\node[scale=0.9] at(0.25,-2.2){Rx angle $\alpha$ in °};
		\node[scale=0.9,rotate=90] at(-4.3,0.25){BRCS (dB)};
		\begin{scope}[xshift=0.11cm,yshift=0.05cm]
			\draw[line width=0.35pt,gray,fill=white] (-0.17,0.25)rectangle(3,-1.5);
			\def\d{0.33}
			\definecolor{RISsim}{rgb}{1,0.5,0.5}
			\definecolor{RISapprox}{rgb}{1,0.1,0.1}
			\definecolor{RISmeas}{rgb}{0.9,0.0,0.0}
			\def\legendscale{0.77}
			\draw[RISsim,densely dashed,line width=0.95pt] (0,0)--+(0.25,0)
				node[right,text=black,scale=\legendscale]{RIS Simulation\vphantom{Ap}};
			\draw[RISapprox,line width=0.95pt] (0,-\d)--+(0.25,0)
				node[right,text=black,scale=\legendscale]{RIS Approximation};
			\draw[RISmeas,line width=1.7pt] (0,-2*\d)--+(0.25,0)
				node[right,text=black,scale=\legendscale]{RIS Measurement\vphantom{Ap}};
			\draw[blue,densely dashed,line width=0.95pt] (0,-3*\d)--+(0.25,0)node[right,text=black,scale=\legendscale]{Reflector Simulation\vphantom{Ap}};
			\draw[blue,line width=1.7pt] (0,-4*\d)--+(0.25,0)node[right,text=black,scale=\legendscale]{Reflector Measurement\vphantom{Ap}};
		\end{scope}
		\end{tikzpicture}}
		\vspace*{-0.5mm}
	\caption{Practical setup and results: (a) shows the realized $7\times2$-element RIS with electronics, (b) depicts the entire measurement setup with the Tx/Rx ridged horns and (c) presents a comparison of the resulting BRCS from simulation, approximation and measurement for a Tx at $\beta=30$° and the RIS optimized for an Rx at $\alpha=0$°. As reference, the BRCS of a plane copper reflector of the same dimensions is included, revealing that even such a small RIS can have a significant impact, exceeding $+15$\,dB @ $\alpha=0$°.}
\end{figure}

\vspace*{4pt}
\section{Summary \& Conclusion}

A simple far-field approximation is introduced, capable of efficiently calculating the full scatter matrices for arbitrary receiver and transmitter positions from only a single full-wave simulation, while maintaining high accuracy suitable for optimization purposes. The validity of this approach is demonstrated through comparisons with optimized capacitance values and is further corroborated by empirical measurements.

\vspace*{4pt}
\section*{Acknowledgment}
This work was carried out as part of the collaborative CHIST-ERA project ”Towards Sustainable ICT: SUNRISE”. 
The authors of OST are supported by the Swiss National Science Foundation (SNF) with Grant 203784.


\vspace*{4pt}
\bibliographystyle{myIEEEtran}
\bibliography{IEEEabrv.bib,aps_2024_ris2.bbl}

\end{document}